\documentclass[11pt]{amsart}

\newcommand{\OO}{\mathcal O}
\newcommand{\Hom}{\operatorname {Hom}}

\newcommand{\Zp}{{\mathbf Z}_{p}}
\newcommand{\Qp}{{\mathbf Q}_{p}}
\newcommand{\FF}{{\mathbf F}_{q}}
\newcommand{\F}{{\mathbf F}}
\newcommand{\Fp}{{\mathbf F}_p}
\newcommand{\CC}{\mathbf C}
\newcommand{\RR}{\mathbf R}
\newcommand{\ZZ}{\mathbf Z}
\newcommand{\QQ}{\mathbf Q}
\newcommand{\mm}{{\mathfrak m}}
\newcommand{\TT}{T}
\newcommand{\UU}{T}
\newcommand{\thet}{\underline{\mathcal D}}
\newcommand{\dd}{\underline{d}}
\newcommand{\DD}{{\mathcal D}}
\newcommand{\Omeg}{{\mathcal B}}
\newcommand{\RD}{{\mathbf D}}

\newcommand{\LC}{\operatorname {Hom}}
\newcommand{\NN}{\mathbf N}
\newcommand{\Maps}{\operatorname {Maps}}
\newcommand{\Ker}{\operatorname {Ker}}
\newcommand{\ord}{\operatorname {ord}}
\newtheorem{theorem}{Theorem}
\newtheorem{corollary}{Corollary}
\newtheorem{lemma}{Lemma}
\newcommand{\eqdef}{:=}
\newcommand{\Teich}{\operatorname {Teich}}

\newcommand{\res}{\operatorname {res}}
\newcommand{\Sym}{\operatorname {Sym}}
\newcommand{\Atop}[2]{\genfrac{}{}{0pt}{}{#1}{#2}}

\title{The Digit Principle}
\author{Keith Conrad}
\address{Department of Mathematics, 
Ohio State University, Columbus, OH 43210-1174}
\email{kconrad@math.ohio-state.edu}

\begin{document}

\begin{abstract}
A number of constructions in function 
field arithmetic involve extensions from linear 
objects using digit expansions.  
This technique is described here as a method of constructing 
orthonormal bases in spaces of continuous functions.
We illustrate 
several examples of orthonormal bases from this viewpoint, 
and we also obtain 
a concrete model for the  
continuous functions on 
the integers of a local field as a quotient of a Tate 
algebra in countably many variables.

\end{abstract}

\maketitle

\section{Introduction}

In function field arithmetic, there is a standard 
construction in which linear objects are extended by using 
digit expansions.  

The Carlitz polynomials are a basic 
example.  Let $\F_{r}[\TT]$ be the polynomial ring in $\TT$ over 
the finite field $\F_{r}$, $\F_{r}[\TT]^{+}$ the subset of monic polynomials. 
For an integer $j \geq 0$, set 
$$
e_{j}(x) \eqdef 
\prod_{\Atop{h \in \F_{r}[\TT]}{\deg(h) < j}} (x - h) \in \F_{r}[\TT][x], 
\ \ \  D_j \eqdef \prod_{\Atop{h \in \F_{r}[\TT]^{+}}
{\deg(h) = j}} h \in \F_{r}[\TT], \ \ \ 
E_{j}(x) \eqdef \frac{e_j(x)}{D_j}.
$$

The polynomial $h = 0$ is included in the product 
defining $e_j(x)$ when $j > 0$, 
and $e_0(x) = 1$.
Since $E_1(x) = 
(x^r-x)/(\TT^r-\TT)$ and 
$h(\TT)^r - h(\TT)$ has all of $\F_{r}$ as roots for any $h(T) \in \F_{r}[T]$, 
$E_1(h(\TT)) \in \F_{r}[\TT]$.  More generally, $e_j(x)$ and 
$E_j(x)$ are both $\F_{r}$-linear maps from $\F_{r}[T]$ to $\F_{r}[T]$.
For any monic $h$ of degree $j$, $e_j(h) = D_j$, so $E_j(h) = 1$.
If $\deg(h) < j$, then $e_j(h) = E_j(h) = 0$.

From the sequences $\{e_j(x)\}$ and $\{E_j(x)\}$ of $\F_{r}$-linear
functions, the Carlitz polynomials are constructed as 
$$
G_{i}(x) \eqdef \prod_{j=0}^{k} e_{j}(x)^{c_j},  \ \ \ \ 
\mathcal{E}_{i}(x) \eqdef \prod_{j=0}^{k} E_{j}(x)^{c_j} = 
\prod_{j=0}^{k} \left(\frac{e_j(x)}{D_j}\right)^{c_j},
$$
where $i = c_0 + c_1r + \dots + c_kr^k$, $0 \leq c_j \leq r-1$.  
Note $E_j(x) = \mathcal{E}_{r^j}(x)$.  The denominator 
$\prod_{j=0}^{k}D_j^{c_j}$ of ${\mathcal E}_i(x)$ is the Carlitz factorial 
${\it \Pi}(i)$.  Basic properties of the Carlitz functions 
can be found in Goss \cite{goss}, \cite[Chap. 3]{gossbk}.

An important property of the Carlitz 
polynomials $\mathcal{E}_{i}(x)$ is 
that every continuous function 
$f \colon \F_{r}[[\TT]] \rightarrow \F_{r}((\TT))$ can be written uniquely 
in the form
\begin{equation}\label{contexp}
f(x) = \sum_{i \geq 0} a_{i}\mathcal{E}_{i}(x), 
\end{equation}
where $a_i \in \F_{r}((\TT))$ and $a_i \rightarrow 0$ as 
$i \rightarrow \infty$.  This is due to Wagner \cite{wagner}, 
who shows as a corollary that 
every continuous $\F_{r}$-linear function $g \colon \F_{r}[[\TT]] \rightarrow 
\F_{r}((\TT))$ can be written uniquely in the form 
\begin{equation}\label{lincontexp}
g(x) = 
\sum_{j \geq 0} c_j{\mathcal E}_{r^j}(x) =  \sum_{j \geq 0} c_{j}E_{j}(x), 
\end{equation}
where $c_j \rightarrow 0$.  The theme which will be seen in 
several guises in this paper is that 
in such situations it 
is simpler to verify an expansion property like (\ref{lincontexp}) 
for linear 
continuous functions first. 
An expansion property like (\ref{contexp}) for general continuous functions 
then follows by an argument that involves little which is special 
about the Carlitz functions ${\mathcal E}_i$ except for their construction 
from digit expansions.  This applies to 
several examples besides the Carlitz basis.  One of these 
examples, due to Baker, yields an interesting 
model for the algebra of continuous functions on the integers of a local 
field.

In characteristic 0, Mahler's theorem says that the binomial 
polynomials $\binom{x}{n} \in \QQ[x]$ are a basis for the 
continuous functions from $\Zp$ to $\Qp$ for {\it all} primes $p$. 
We will consider some analogues of this phenomenon in 
positive characteristic, requiring a passage at times 
between global fields and their completions.   
Our notational conventions in this regard 
are as follows.
The global fields we will consider in positive 
characteristic will be of the form 
$\F_{r}(T)$, whose 
completion at any place has the form $\FF((u))$ 
for some finite field $\FF$ and uniformizing parameter 
$u$. 
We also write the residue field as $\F_{u}$.

We need some additional notation for spaces of maps.  

For any local field $K$ (always nonarchimedean) we denote its 
integer ring and the corresponding maximal ideal as 
$\OO$ and $\mm$.  The residue field is denoted $\F$.
We write $C(\OO,K)$ for the continuous functions from $\OO$ to 
$K$, topologized with the sup-norm.  We similarly define 
$C(\OO,\OO)$ and $C(\OO,\F)$, viewing $\F$ as 
a discrete space.  So any 
element of $C(\OO,\F)$ factors through a finite quotient of $\OO$.

When $K$ is a local field of positive characteristic, 
we write 
$\Hom_{\F}(\OO,K), \Hom_{\F}(\OO,\OO)$, and $\Hom_{\F}(\OO,\F)$ 
for the {\it continuous} $\F$-linear maps from $\OO$ 
to the corresponding sets.  (In particular, 
continuous $\F$-linear maps from $\OO$ to $\F$ always factor 
through some finite quotient.)  We will at times consider 
linear maps relative to a subfield $\F' \subset \F$, 
so write $\Hom_{\F'}$ in these cases. Note $\Hom_{\F}$ 
and $\Hom_{\F'}$ will 
never mean algebra homomorphisms.

For finite sets $A$ and $B$, $\Maps(A,B)$ is the set of all 
functions from $A$ to $B$.  This will only arise when 
$B$ is a finite field, making the space of maps an 
vector space in the natural way.

I thank D. Goss and W. Sinnott for discussions on the topics in this 
paper.

\section{Background}

Let $(E,||\cdot||)$ a Banach space over a local field $K$. 
Let 
$$
E_0 \eqdef \{x \in E : ||x|| \leq 1\}.
$$
So, using the notation given in the introduction, the {\it residual space}  
$\overline{E} \eqdef E_0/\mm E_0$ is a vector space over 
the residue field $\F = \OO/\mm$. 

We assume throughout that every 
nonzero element of $E$ has its norm value in the value 
group of $K$.  This is required in order 
to know that all elements of 
$E$ can be scaled to have norm 1, and in particular 
that $\mm E_0 = \{x \in E : ||x|| < 1\}$.

\hfill

\noindent
{\bf Example.}
Let 
$C(\OO,K)$ be the $K$-Banach space of continuous functions 
from $\OO$ to $K$, topologized by the sup-norm.  
Since we use the sup-norm, the space 
$E = C(\OO,K)$ has $||E|| = |K|$ and 
$$
\overline{C(\OO,K)} \cong C(\OO,\F).
$$

\hfill

\noindent
{\bf Example.}
Let $K$ have positive characteristic, so the residue field 
$\F$ is a subfield of $\OO$.  We consider $E = \LC_{\F}(\OO,K)
\subset C(\OO,K)$. 
Again $||E|| = |K|$.  Since $\F \subset \OO$, 
$$
\overline{\LC_{\F}(\OO,K)} 
\cong \LC_{\F}(\OO,\F).
$$

\hfill

A sequence $\{e_0, e_1, e_2, \dots\}$ in $E$ is 
called an {\it orthonormal 
basis} if each $x \in E$ has a representation as 
$$
x = \sum_{n \geq 0} c_{n}e_{n},
$$
where $c_n \in K$ with $c_n \rightarrow 0$ and
$$
||x|| = \max_{n \geq 0} |c_n|.
$$
The coefficients in such a representation are unique.

\begin{lemma}\label{orc}
For a local field $(K,|\cdot|)$ and a $K$-Banach space $(E,||\cdot||)$, 
where $||E|| = |K|$, 
a necessary and sufficient 
condition for a sequence $\{e_n\}$ in $E$ to be an orthonormal 
basis is that 
every $e_n$ lies in $E_0$ and the reductions 
$\overline{e}_n \in \overline{E}$ 
form 
an $\F$-basis of $\overline{E}$ in the algebraic sense, 
i.e., using finite linear combinations.
\end{lemma}

\begin{proof} 
See Serre \cite[Lemme I]{serre} or 
Lang \cite[\S 15.5]{lang}.
\end{proof}

For a counterexample to Lemma \ref{orc}
when $K$ is a non-archimedean complete field with a 
non-discrete valuation, see 
Bosch, G\"untzer, Remmert \cite[p. 118]{bgr} or 
van Rooij \cite[p. 183]{vr}

Our use of the notation 
$e_n$ for a vector in a Banach space 
should not be confused with the Carlitz polynomial 
written as $e_n(x)$.  We will only use the 
Carlitz polynomial $e_n(x)$ again within the proofs of 
Lemma \ref{cbasis} and Lemma \ref{prep}.

By Lemma \ref{orc}, 
functions $e_i$ in $C(\OO,K)$ form an orthonormal basis if and only if 
they map $\OO$ to $\OO$ and their reductions $\overline{e}_i =  
e_i \bmod \mm$ 
are an algebraic basis of 
$$
C(\OO,\F) = \varinjlim \Maps(\OO/\mm^n,\F).
$$
So the construction of an orthonormal 
basis of $C(\OO,K)$ is reduced to a linear algebra problem: 
verifying a sequence in $C(\OO,\F)$ is an $\F$-basis.  
For example, 
let $q = \#\F$ and suppose 
for all $n \geq 0$ (or simply infinitely many $n \geq 0$) 
that the functions 
$\overline{e}_0,\dots,\overline{e}_{q^n-1} \colon \OO \rightarrow \F$ 
are well-defined modulo $\mm^n$ and give 
an $\F$-basis of $\Maps(\OO/\mm^n,\F)$.  Then 
the set of all $e_i$ forms an orthonormal basis of $C(\OO,K)$.  
We will often intend this particular remark when we refer later 
to Lemma \ref{orc}.

The standard 
examples, such as the binomial 
polynomials $\binom{x}{n}$ viewed in $C(\Zp,\Fp)$ and 
the Carlitz polynomials ${\mathcal E}_{n}(x)$ 
viewed in 
$C(\FF[[T]],\FF)$, are usually checked to be 
bases by combinatorial inversion 
formulas involving certain sequences of difference operators.
We will not utilize any difference operators, although 
implicitly they provide one way of checking the 
binomial and Carlitz polynomials take integral values.

The case we are interested in first 
is local fields of positive characteristic.  
Let $K$ be such a local field, with $\OO$ its ring of 
integers and $\F$ the residue field.
Rather than starting with $C(\OO,K)$, we begin 
with the closed subspace $\LC_{\F}(\OO,K)$ of continuous 
$\F$-linear functions from $\OO$ to $K$.

A sequence $e_j$ in $\LC_{\F}(\OO,K)$ consisting 
of functions sending $\OO$ to $\OO$ is an 
orthonormal basis of $\LC_{\F}(\OO,K)$ if and only if
the reductions 
$\overline{e}_j$ form an  
algebraic basis of the residual space $\LC_{\F}(\OO,\F)$.
Let $e_0,e_1,e_2,\dots$ be an orthonormal basis of 
$\LC_{\F}(\OO,K)$, and $q = \#\F$.  We define 
a sequence of functions $f_i$ for $i \geq 0$ by writing $i$ in base $q$ as 
$$
i = c_0 + c_1q + \dots + c_{n-1}q^{n-1}, \ \ 0 \leq c_j \leq q-1,
$$
and then setting 
\begin{equation}\label{fecon}
f_i \eqdef e_{0}^{c_0}e_{1}^{c_1}\dots e_{n-1}^{c_{n-1}}.
\end{equation}
Note $e_j = f_{q^j}$.  If $c = 0$, $e_j^c$ is the function 
that is identically 1, even if $e_j$ vanishes somewhere.  
The construction in (\ref{fecon}) will be called the 
extension of the $e_j$ by digit expansions, or the 
extension by $q$-digits 
if the reference to $q$ is worth clarifying.

We show in the next section that the 
$f_i$ form an orthonormal basis of $C(\OO,K)$, a fact which 
we refer to as the ``digit principle.''

\section{Extending an Orthonormal Basis}

\begin{theorem}[Digit Principle in Characteristic $p$]\label{exdigit}
Let $K$ be a local field of positive characteristic, with 
integer ring $\OO$ and residue field $\F$ of size $q$.
The extension of an orthonormal basis of $\LC_{\F}(\OO,K)$ via 
$q$-digit 
expansions 
produces an 
orthonormal basis for $C(\OO,K)$.
\end{theorem}

\begin{proof}
Let $\{e_j\}_{j \geq 0}$ be an orthonormal basis of 
$\LC_{\F}(\OO,K)$, so $\{\overline{e}_j\}_{j \geq 0}$ is 
an $\F$-basis of 
$$
\overline{\LC_{\F}(\OO,K)} = \LC_{\F}(\OO,\F) = 
\bigoplus_{j \geq 0} \F\overline{e}_j.
$$
Let $H_n = \cap_{j=0}^{n-1} \Ker(\overline{e}_j)$, so 
$H_n$ is a closed subspace of $\F$-codimension $n$ in $\OO$, 
$H_{n+1} \subset H_n$, and 
$\cap H_n = 0$.  Therefore $\OO \cong \varprojlim \OO/H_n$, so 
$C(\OO,\F) = \varinjlim \Maps(\OO/H_n,\F)$.   Viewing 
$\overline{e}_0,\dots,\overline{e}_{n-1}$ as functions on 
$\OO/H_n$, they form an 
$\F$-basis of the $\F$-dual space $(\OO/H_n)^{*}$.
So we are reduced to an 
issue about linear algebra over finite fields: for 
$q = \#\F$ and $0 \leq i \leq q^n-1$, 
do the $q^n$ 
reduced functions $\overline{f}_i$, as constructed in (\ref{fecon}), 
form a basis of $\Maps(\OO/H_n,\F)$?

Let $V$ be a finite-dimensional 
$\FF$-vector space, of dimension (say) $n$.  Let 
$\varphi_0,\dots,\varphi_{n-1}$ be a basis of $V^{*}$.
Extend the $\varphi_j$ to a set of 
polynomial functions on $V$ by using digit 
expansions.  That is, 
for $0 \leq i \leq q^n-1$ write $i = c_0 + c_1q + \dots + c_{n-1}q^{n-1}$ 
in base $q$ and set
$$
\Phi_i = \varphi_{0}^{c_0}\cdot \dots \varphi_{n-1}^{c_{n-1}}.
$$
So $\varphi_{j} = \Phi_{q^j}$ and $\Phi_{0} = 1$.
By a dimension count, we just need to 
show the functions $\Phi_i$ are a basis of $\Maps(V,\FF)$.
It suffices to show the $\Phi_i$ span $\Maps(V,\FF)$.

Let 
$$
\{v_0,v_1,\dots,v_{n-1}\} \subset V
$$ 
be the dual basis 
to the $\varphi_j$.  For $v \in V$, write
$$
v = a_0v_0 + a_1v_1 + \dots + a_{n-1}v_{n-1}, 
$$
where $a_j \in \FF$. 
Taking an idea from 
the proof of the Chevalley-Warning 
Theorem in Serre \cite[p. 5]{cia}, define $h_v \colon V \rightarrow \FF$ by 
$$
h_v(w) \eqdef \prod_{j=0}^{n-1} (1 - (\varphi_j(w) - a_j)^{q-1}) = 
\prod_{j=0}^{n-1}(1 - (\varphi_j(w) - \varphi_j(v))^{q-1}).
$$
Since $h_v(w)$  is 1 when $w = v$ and $h_{v}(w) = 0$ when $w \not= v$, 
the $\FF$-span of the $h_v$ is all of  $\Maps(V,\FF)$.
Expanding the product defining $h_v$ shows 
$h_v$ is in the span of the $\Phi_i$ since the exponents of 
the $\varphi_j$ in the product never exceed $q-1$.
This concludes the proof.
\end{proof}

In terms of a basis $\varphi_0,\dots,\varphi_{n-1}$ of $V^{*}$, 
the main point of the 
proof is that the natural map
\begin{equation}\label{symis}
\FF[\varphi_0,\dots,\varphi_{n-1}]/(\varphi_0^{q}-\varphi_0,\dots,
\varphi_{n-1}^{q}-\varphi_{n-1}) \rightarrow \Maps(V,\FF)
\end{equation}
is an isomorphism.  This is the familiar fact that any 
function $\FF^n \rightarrow \FF$
has the same graph as 
a polynomial whose variables all have degree at most $q-1$.  
Without reference to a basis of $V^{*}$, (\ref{symis}) 
can be written as $\Sym(V^{*})/I \cong \Maps(V,\FF)$, 
where $I$ is the ideal generated by all $g^q - g$, $g \in \Sym(V^{*})$.

In practice, the codimension condition at the start of 
the proof of Theorem \ref{exdigit}
may be known not because the $\overline{e}_j$ are an $\F$-basis 
of $\Hom_{\F}(\OO,\F)$
but by some other means in the course of showing 
these functions are a basis.

Although Theorem \ref{exdigit}
gives a natural explanation for some aspects of 
digit expansions in function field arithmetic, there are settings 
where the use of digit expansions remains mysterious.  For instance, 
is there a 
natural explanation for the role of digit expansions in 
the construction of function 
field Gamma functions (cf. Goss \cite{goss}, \cite{gossbk})?

In the notation of Theorem \ref{exdigit}, 
let $\F'$ be a subfield of $\F$ and 
consider the closed subspace $\LC_{\F'}(\OO,K)$ 
of $\F'$-linear continuous functions from $\OO$ to $K$.
Since $\F \subset \OO$, 
the residual space 
$\overline{\LC_{\F'}(\OO,K)}$ equals $\LC_{\F'}(\OO,\F)$. 
For any $e \in \LC_{\F'}(\OO,\OO)$, note the 
kernel of $\overline{e} \colon \OO \rightarrow \F$ is 
not typically an $\F$-subspace, only an $\F'$-subspace.  
By imposing a condition on 
the kernel of $\overline{e}$ which is 
automatically satisfied when $\F' = \F$, we can extend 
the scope of Theorem \ref{exdigit} as follows.

\begin{theorem}\label{rrq}
Let $K$ be a local field of positive characteristic, with integer 
ring $\OO$ and residue field $\F$.
Let $\F'$ be a subfield of $\F$, with $\F' = r$ and $\F = q = r^d$.
If $\{e_j\}_{j \geq 0}$ is an orthonormal basis of $\LC_{\F_r}(\OO,K)$ 
such that $\cap_{j=0}^{dn-1} \Ker(\overline{e}_j)$ has $\F'$-codimension 
$dn$ in $\OO$ for all $n \geq 1$, 
then the extension of the $e_j$ by $r$-digits gives an orthonormal 
basis of $C(\OO,K)$.
\end{theorem}

Note the digit extension in the theorem is by $r$-digits, not by $q$-digits 
(where $q = r^d$).

\begin{proof}
Let $H_n = \cap_{j=0}^{dn-1} \Ker(\overline{e}_j)$, so 
$\#(\OO/H_n) = r^{dn} = q^n$. 
For $0 \leq j \leq dn-1$, the maps 
$\overline{e}_j \colon \OO  \rightarrow \F$ 
give well-defined 
$\F'$-linear maps from $\OO/H_n$ to $\F$.
By hypothesis they are linearly independent over $\F$, 
and since $\Hom_{\F'}(\OO/H_n,\F)$ has dimension $dn$ as an 
$\F$-vector space (indeed, 
$$
\dim_{\F}(\Hom_{\F'}(\OO/H_n,\F)) = \dim_{\F'}(\Hom_{\F'}(\OO/H_n,\F')) 
= \dim_{\F'}(\OO/H_n) = dn),  
$$
the functions $\overline{e}_0,\dots,\overline{e}_{dn-1}$, 
when viewed in $\Hom_{\F'}(\OO/H_n,\F)$, are 
an $\F$-basis.  Therefore 
the functions $\overline{e}_0,\dots,\overline{e}_{n-1}$
separate the points of $\OO/H_n$.  
(Intuitively, this 
situation as analogous to 
a finite-dimensional $\RR$-vector space 
$W$ and a $\CC$-basis $f_1,\dots,f_m$ of $\Hom_{\RR}(W,\CC)$. 
Such a $\CC$-basis separates any two 
points of $W$, since an $\RR$-dual vector $h \colon W \rightarrow 
\RR$ does the job and we view $\RR \subset \CC$ 
to realize $h$ as a $\CC$-linear combination 
of the $f_k$; thus one of the $f_k$ separates the two points.)

An argument as in the proof 
of Theorem \ref{exdigit} 
then shows that $\Maps(\OO/H_n,\F)$ is spanned over $\F$ 
by the monomials
$$
\overline{e}_{0}^{b_0}\cdots 
\overline{e}_{dn-1}^{b_{dn-1}},  \ \ \ \ 0 \leq b_j \leq q-1.
$$
This set has size $q^{dn}$, which is too large 
(when $d > 1$) to be an $\F$-basis 
of $\Maps(\OO/H_n,\F)$. 

To cut down the size of this spanning set, note 
any 
$e_{j}^{r^k}$ is $\F'$-linear, so 
in $\Maps(\OO/H_n,\F)$ we can write 
$\overline{e}_j^{r^k}$ 
as an 
$\F$-linear combination of $\overline{e}_0,\dots,\overline{e}_{dn-1}$.
Therefore for all $n \geq 1$, $\Maps(\OO/H_n,\F)$ is spanned over $\F$ by 
$$
\overline{e}_{0}^{c_0}\cdots 
\overline{e}_{dn-1}^{c_{dn-1}},  \ \ \ \ 0 \leq c_j \leq r-1, 
$$
so this set is an $\F$-basis.  We're done by Lemma \ref{orc}.
\end{proof}

As formulated so far, the digit principle does not apply 
in characteristic 0 since there is no analogue in characteristic 
0 of the subspace of linear functions.  
However, a remark of Baker \cite[p. 417]{baker} shows that 
replacing the 
linear condition with a property that comes up in the proof of 
Theorem \ref{exdigit} extends 
the digit principle to characteristic 0, as follows.

\begin{theorem}[Digit Principle in any Characteristic]\label{char0}
Let $K$ be any local field, 
$\OO$ its ring of 
integers, $\F$ the residue field, and $q = \#\F$.  
Let $H_n$ be a sequence of open subgroups of $\OO$ such that 
$H_{n+1} \subset H_n$ and $\cap H_n = 0$.
Suppose there is a 
sequence $e_0,e_1,e_2,\dots$ in $C(\OO,K)$ such that each 
$e_j$ maps $\OO$ to $\OO$ and for all $n \geq 1$, 
the 
reductions $\overline{e}_0, \dots, \overline{e}_{n-1} \in C(\OO,\F)$ 
are constant on cosets of $H_n$ and the map 
$$
\OO/H_n \rightarrow \F^n \ \ 
\text{ given by } \ \ x \mapsto (\overline{e}_0(x),\dots,
\overline{e}_{n-1}(x))
$$
is bijective.  
$($So $\#\OO/H_n = q^n$.$)$ 
The extension of the sequence $e_j$ by $q$-digits 
gives an orthonormal basis of 
$C(\OO,K)$.
\end{theorem}

In positive characteristic $H_n$ is an $\F$-vector space, 
so a natural (though not essential) 
way for the functions
$\overline{e}_0,\dots,\overline{e}_{n-1}$ 
on $\OO/H_n$ to satisfy the bijectivity hypothesis is for 
them to be an $\F$-basis of the dual space $(\OO/H_n)^{*}$, which is how 
Theorem \ref{exdigit} is proved.

\begin{proof}
For $v \in \OO/H_n$, 
our hypotheses make the function $h_v \colon \OO/H_n \rightarrow \F$ given by 
$$
h_v(w) = \prod_{j=0}^{n-1} (1 - (\overline{e}_j(w) - \overline{e}_j(v))^{q-1})
$$
equal to 1 for $w = v$ and 0 for $w \not= v$, so the proof 
of Theorem \ref{exdigit} still works.
\end{proof}

Although Theorem \ref{char0} includes the previous theorems 
as special cases, from the viewpoint of applications 
the linearity hypotheses in the earlier theorems make 
it convenient to isolate them separately and independently 
from Theorem \ref{char0}.
(This is partly why they were treated first.)
It might be worth referring to Theorems \ref{exdigit} and \ref{rrq}
as the linear digit principle to distinguish them from Theorem \ref{char0}, 
but we won't adopt this extra appellation here.

While we formulated Theorem \ref{char0} for 
general classes of subgroups $H_n$, in the applications in 
Section \ref{apps} we will 
only encounter $H_n = \mm^n$.

In Theorem \ref{char0} we can consider instead 
a sequence $e_j$ in $C(Z,K)$, where 
$Z = \varprojlim Z_n$ is profinite.
Suppose for all $n \geq 1$ 
that the reduced functions $\overline{e}_j \colon Z \rightarrow \F$ 
for $0 \leq j \leq n-1$ 
are constant on the fibers of the natural projection $Z \rightarrow Z_n$
and the induced map $Z_n \rightarrow \F^n$ given by 
$$
z \mapsto (\overline{e}_0(z),\dots, \overline{e}_{n-1}(z))
$$
is a bijection.  Then the extension of the $e_j$ by $q$-digits 
is an orthonormal basis of $C(Z,K)$.

\section{Hyperdifferential Operators}

Some of the applications we give in Section \ref{apps} involve 
a set of differential operators whose main features 
we summarize here.

For any field $F$ and 
integer $j \geq 0$, 
the $j$th hyperdifferential operator $\DD_j = \DD_{j,T}$, also called the 
divided power   
derivative, acts 
on the polynomial ring 
$F[\UU]$ by $\DD_{j}(\UU^m) = \binom{m}{j}\UU^{m-j}$ 
for $m \geq 0$ and is extended by $F$-linearity to all polynomials.
These operators were first studied by Hasse and Schmidt \cite{hasse}
and Teichm\"uller \cite{teichm}.

If $F$ has characteristic 0 then 
$$
\DD_j = \frac{1}{j!}\frac{d^{j}}{d\UU^{j}}, 
$$
but this formula holds in characteristic $p$ only for 
$j \leq p-1$.
Unlike the ordinary higher derivatives, $\DD_j$ 
is not a trivial operator in characteristic $p$ when 
$j \geq p$.  For example,
$$
\DD_3(1 + \UU + 2\UU^3 + 2\UU^7 + \UU^9) = 
2+ 70\UU^4 + 84\UU^6 \equiv 2 + \UU^4 \bmod 3.
$$
Note the constant term of $\DD_j(f(\UU))$ is simply the $j$th 
Taylor coefficient of $f(\UU)$.  
While $\DD_j$ is not an iterate of $\DD_1$, in characteristic 
$p$ it does share the property that the $p$-fold iterate of 
$\DD_j$ is identically 0.  We will generally not be considering 
iterates of the $\DD_j$, but rather their products in the sense of 
functions.

\begin{theorem}
\label{hyperchar}
The hyperdifferential operators $\DD_j \colon F[\UU] \rightarrow F[\UU]$ 
are the unique sequence of maps such that 

i$)$ each $\DD_j$ is $F$-linear,

ii$)$ $\DD_0\UU = \UU$, $\DD_{1}\UU = 1$, 
and $\DD_{j}\UU = 0$ for $j \geq 2$, 

iii$)$ $($Leibniz Rule$)$ 
For all $j \geq 0$, 
$\DD_j(fg) = \sum_{k=0}^{j} (\DD_{k}f)(\DD_{j-k}g)$ for 
all $f, g$ in $F[\UU]$.
\end{theorem}

\begin{proof}
To check that the hyperdifferential operators 
satisfy these three properties, 
only the third has some (slight) content. 
By $F$-linearity it reduces to 
the case $f = \UU^a$ and $g = \UU^b$, in which case 
the Leibniz rule becomes the Vandermonde formula
$$
\binom{a+b}{j}T^{a+b} = \sum_{k=0}^{j}\binom{a}{k}\binom{b}{j-k}T^{a+b}.
$$

Conversely, properties ii and iii suffice to recover the formula 
$\DD_j(\UU^m) = \binom{m}{j}\UU^{m-j}$ for $m \geq 0$, 
which by property i forces $\DD_j$ to be the 
$j$th hyperdifferential operator.  

\end{proof}

The Leibniz rule extends to more than two factors, as 

\begin{equation}\label{mprod}
\DD_{j}(f_1\cdots f_m) = \sum_{\Atop{k_1 + \dots + k_m = j}{k_1,\dots,
k_m \geq 0}} \DD_{k_1}(f_1)\cdots \DD_{k_m}(f_m).
\end{equation}

Since $\DD_{j} = (1/j!)(d/d\UU)^j$ in characteristic 0, 
it is natural over any $F$ to view 
the operators $\DD_j$ as coefficient functions of a formal 
Taylor expansion
\begin{equation}\label{texp}
\thet \colon 
f \mapsto \sum_{j \geq 0} (\DD_{j}f)X^j
\end{equation}
from $F[\UU]$ to $F[\UU][[X]]$.  (The image of $\thet$ here is only in 
$F[\UU][X]$, but it is convenient for later extension problems 
to have the target space be formal 
power series in $X$.)
Properties i and iii of Theorem \ref{hyperchar} are equivalent to 
$\thet$ being 
an $F$-algebra homomorphism.  
Property ii just says $\thet(\UU) = \UU + X$.
By a computation, 
$\thet(\UU^m) = (\UU+X)^m$, so $\thet$ is indeed an $F$-algebra 
homomorphism, in fact it is simply 
given by $\thet(f(\UU)) = f(\UU + X)$. 
This is an alternate proof of Theorem 
\ref{hyperchar}.

For any $f, g \in F[\UU]$, consider 
the expression for the coefficient of $X^j$ in  
$$
\thet(f^ng) = \thet(f)^n\thet(g)
$$ 
as a sum of monomials arising from multiplication of the series 
on 
the right side.   To obtain a coefficient of $\UU^j$ 
by selecting one term from each of these series, 
at least $n-j$ of the $n$ factors equal to $\thet(f)$ must contribute 
their constant term, which is $f$.  Therefore 
\begin{equation}\label{hypercong2}
\DD_{j}(f^ng) \equiv 0 \bmod f^{n-j}
\end{equation}
for $n \geq j$.  
In particular, each $\DD_{j}$ is $f$-adically continuous 
for any polynomial $f$ in $F[\UU]$.
Since 
$\thet(f) \equiv f(\UU) + f'(\UU)X \bmod X^2$, we have
$\thet(f) \equiv f'(\UU)X \bmod (X^2,f)$, 
so looking at the coefficient of $X^j$ in $\thet(f^j) = \thet(f)^j$ shows
\begin{equation}\label{hypercong1}
\DD_{j}(f^j) \equiv (f')^j \bmod f.
\end{equation}

The sequence of operators $\DD_j$ on $F[T]$ is a special case of a higher 
derivation, which we now recall. 

For any commutative ring $R$ and $R$-algebra $A$, 
a {\it higher $R$-derivation} on $A$ is a sequence of 
$R$-linear maps $d_j \colon A \rightarrow A$ for $j \geq 0$ such that 
$d_0$ is the identity map and $d_j(ab) = \sum_{k=0}^{j} d_k(a)d_{j-k}(b)$ 
for all $a, b \in A$.  (So $d_1$ is a derivation in the usual sense.)
Equivalently, the map $\dd \colon A \rightarrow A[[X]]$ 
given by $\dd(a) =  \sum_{j \geq 0} d_{j}(a)X^j$ 
is an $R$-algebra homomorphism that is a section to the 
canonical map $A[[X]] \rightarrow A$.  This equivalent viewpoint shows  
that any higher $R$-derivation $\{d_j\}$ on $A$ extends uniquely 
to a higher $R$-derivation on any localization $S^{-1}A$ of $A$; 
we simply extend the corresponding $R$-algebra map $\dd$ 
uniquely to an $R$-algebra map from $S^{-1}A$ to $(S^{-1}A)[[X]]$.
For example, the sequence of hyperdifferential operators 
$\DD_j$ on $F[\UU]$ 
extends uniquely to a higher $F$-derivation on 
$F(\UU)$.  For nonzero $f$ in $F[\UU]$, 
the Leibniz rule computes a formula for $\DD_{j}(1/f)$ 
inductively, though using such formulas to prove the $\DD_j$ 
satisfy the Leibniz rule on the field $F(\UU)$ would be a terrific mess.
(Remember that the $\DD_j$ are not iterates of $\DD_1$.)

\begin{theorem}
Let $K$ be a field, $F$ a subfield.  Any higher $F$-derivation 
on $K$ extends uniquely to a higher $F$-derivation on any 
separable algebraic extension $L/K$.
\end{theorem}

In particular, the only higher $F$-derivation $\{d_n\}$ on a separable 
algebraic extension of $F$ is given by $d_n = 0$ for $n \geq 1$.

\begin{proof}
It suffices to assume $L/K$ is a finite extension, say 
$L = K(\alpha_0)$ where $\alpha_0$ is the root of the 
separable monic irreducible polynomial $\pi(Y) \in K[Y]$.

Let $\underline{\delta} 
\colon K \rightarrow K[[X]]$ be a higher $F$-derivation 
on $K$.  Any extension of $\underline{\delta}$ to a higher $F$-derivation 
on $L$ must send $\alpha_0$ to an element of $L[[X]]$ 
which has constant term $\alpha_0$ and is a root to the 
polynomial $\pi^{\underline{\delta}}(Y) \in K[[X]][Y]$, 
where $\pi^{\underline{\delta}}(Y)$ is obtained by 
applying $\underline{\delta}$ to the coefficients of $\pi(Y)$.
This polynomial is 
irreducible over $K((X))$ by Gauss' Lemma.
Since $\pi^{\underline{\delta}}(Y) \bmod X = \pi(Y)$ has 
$\alpha_0$ as a simple root in the residue field 
$L[[X]]/(X) = L$, $\alpha_0$ lifts uniquely 
to a root $\alpha(X) \in L[[X]]$ 
by Hensel's Lemma.  So $\underline{\delta}$ extends uniquely 
to an $F$-algebra map $L \rightarrow L[[X]]$ 
by sending $\alpha_0$ to $\alpha(X)$.
\end{proof}

Taking $f = \UU$ in (\ref{hypercong2}), 
all $\DD_j$ extend by $T$-adic continuity to $F[[\UU]]$, 
providing $F[[\UU]]$ with a higher $F$-derivation, which 
then extends uniquely to a higher $F$-derivation on 
$F((\UU))$.  In particular, 
this extension of $\thet$ 
to $F((\UU))$ is given by 
$$
\thet\left(\frac{1}{T}\right) = 
\frac{1}{\UU + X} = \sum_{j \geq 0} \frac{(-1)^{j}}{\UU^{j+1}}X^j, 
$$
which upon raising to the $m$th power shows 
$\DD_j(\UU^m) = \binom{m}{j}\UU^{m-j}$ for all $m \in \ZZ$,

\begin{theorem}\label{Djrcomp}
Let $v$ be any place of $F(\UU)$ which is trivial on $F$ 
and has a residue field which is separable over $F$.
The maps $\DD_{j}$ on 
$F(\UU)$ extend continuously to 
the completion of $F(\UU)$ at $v$, where they form a 
higher $K$-derivation for $K$ any coefficient 
field in the completion such that $K$ contains $F$.
\end{theorem}

We have already seen this for $v$ the $\UU$-adic place.  If 
$F = \FF$ is a finite field, then $v$ can be any place on $\FF(\UU)$, 
and we can canonically take the residue field $\F_v$ of $v$ 
to be the coefficient field of the completion. So the maps $\DD_j$ on 
$\FF(\UU)$ extend by continuity to a higher $\F_v$-derivation 
on the completion at $v$.

\begin{proof}
We take two cases, depending on whether 
$v$ corresponds to a monic separable irreducible polynomial in 
$F[\UU]$ or to $1/\UU$.

If $v$ is a place corresponding to a monic 
separable irreducible $\pi$ in $F[\UU]$, 
(\ref{hypercong2}) shows that the 
$\DD_{j}$ are all $v$-adically continuous, so 
they all extend by continuity to the completion 
$\OO \eqdef \widehat{F[\UU]}_{v}$ 
and still satisfy $F$-linearity and the Leibniz rule.
The corresponding $F$-algebra 
homomorphism $\thet \colon \OO \rightarrow \OO[[X]]$
given by $\thet(g) = \sum (\DD_{j}g)X^j$ is a 
higher $F$-derivation on $\OO$.

Since $\pi$ is separable, 
$\OO$ has a coefficient field, say $K$, which contains $F$. 
Since $K/F$ is separable, 
the restriction of $\thet$ to $K$ must be the usual inclusion 
$K \hookrightarrow K[[X]]$, so $\thet$ is a higher $K$-derivation on 
$\OO$ and therefore also on the fraction field of $\OO$.

If $v$ is the place corresponding to $1/\UU$, we set
$S = 1/\UU$ and note that $\DD_{j}(S^m) = \binom{-m}{j}S^{m+j}$. 
So the $\DD_{j}$ are $S$-adically continuous on 
$F[S] = F[1/\UU]$.  They form a higher $F$-derivation on 
$F[S]$, since this is a subfield of $F(\UU)$.
The continuous extension of all $\DD_{j}$ to the completion 
$F[[1/\UU]]$ (by continuity) and then to 
$F((1/\UU))$ (by algebra) is along similar lines to the previous case.
\end{proof}

As references for additional properties of higher derivations, 
see Okugawa \cite{okug} and 
Kawahara and Yokoyama \cite{kawa}.  Okugawa includes 
an additional condition on a higher derivation 
$\dd \colon A \rightarrow A[[X]]$, namely that 
$d_j \circ d_k = \binom{k+j}{j} d_{j+k}$.
This is motivated by 
the composition rule for hyperdifferential operators on 
$F[T]$.  With this additional condition as part of the definition,  
the above results on extending higher derivations remain 
true, but the proofs involve 
some further calculations.

\section{Examples}
\label{apps}

We now apply the digit principle to compute several examples of 
orthonormal bases on spaces of continuous functions. 
There will be no discussion
of a corresponding difference calculus which gives a formula for 
the coefficients in such a basis.

\begin{lemma}
\label{cbasis}
Let $K = \FF((\TT)), \OO = \FF[[\TT]]$.
The $\FF$-linear Carlitz 
polynomials $E_i(x)$ are an orthonormal basis 
for all $\FF$-linear continuous functions from $\OO$ to $K$.
\end{lemma}

\begin{proof}
By Lemma \ref{orc}, it suffices to show the 
reductions $\overline{E}_j(x)$ are an algebraic basis of 
the space of continuous $\FF$-linear maps from $\OO$ to $\FF$.
We show $\overline{E}_0,\dots,\overline{E}_{n-1}$ form a 
basis of the $\FF$-dual space $(\OO/T^n)^{*}$ for all $n$.

For $0 \leq j < n$, $E_j(\TT^n) \equiv 0 \bmod \TT$ since 
$$
\ord_{\TT}(e_j(\TT^n)) > \ord_\TT(D_j) = 1 + q + q^2 + \dots + q^{j-1}.
$$
Indeed, by the definition of $e_j(x)$, when $n > j$ 
\begin{eqnarray*}
\ord_\TT(e_j(\TT^n)) & = & n + \sum_{k=0}^{j-1}\sum_{\Atop{h \in \FF[\TT]}
{\deg(h) = k}} \ord_{\TT}(h) \\
& = & n + \sum_{k=0}^{j-1} (q-1)\ord_{\TT}(D_d) \\
& = & n + (1 + q + \dots + q^{j-1}) - j \\
& > & \ord_{\TT}(D_j).
\end{eqnarray*}

So for $0 \leq j \leq n-1$, $\overline{E}_j(x)$
is a well-defined function from  
$\FF[\TT]/\TT^n$ to $\FF$.  Since $E_{j}(x)$ vanishes at 
$1,\TT,\dots,\TT^{j-1}$ and $E_{j}(\TT^j) = 1$,
the $n \times n$ 
matrix $(E_{j}(\TT^k))$ is triangular with 
1's along the main diagonal, so it is invertible.
Reducing the matrix entries from $\OO$ into $\FF$ gives an 
invertible matrix, so 
$\overline{E}_0,\dots,\overline{E}_{n-1}$ 
forms a basis of $(\OO/T^n)^{*}$ for all $n$.
\end{proof}

\begin{theorem}
\label{cw1}
Let $K = \FF((\TT)), \OO = \FF[[\TT]]$.
The Carlitz 
functions ${\mathcal E}_i(x)$ are an orthonormal basis of 
the continuous functions from $\OO$ to $K$.
\end{theorem}

\begin{proof}
Use Lemma \ref{cbasis} and the digit principle.
\end{proof}

For a finite field $\F_r$, the construction of the 
Carlitz polynomials
and the hyperdifferential operators on $\F_r[T]$ 
depends on the distinguished generator $T$, and in 
the case of the Carlitz polynomials the construction 
also depends on the 
coefficient field $\F_r$.  To indicate this 
dependence, when it useful, we will write $E_j$ and $\DD_j$ 
as $E_{j,T,r}$ and $\DD_{j,T}$.  This will be necessary when 
we have the Carlitz polynomials or hyperdifferential 
operators that are attached to a global field $\F_r(T)$ act on 
one of the completions $\FF((u))$.  This completion has 
its own local Carlitz polynomials $E_{j,u,q}$ and hyperdifferential 
operators $\DD_{j,u}$ which are typically different from 
the functions $E_{j,T,r}$ and $\DD_{j,T}$ coming from 
the global field.

In Theorem \ref{cw1}, the orthonormal basis 
on $C(\FF[[T]],\FF((T)))$ is constructed via $q$-digits from 
the $\FF$-linear Carlitz polynomials
in $\FF(T)[X]$.  We can instead start with 
a ring $\F_r[T]$, complete it at a prime 
$\pi$, and consider the 
globally constructed $\F_r$-linear 
Carlitz polynomials $E_{j,T,r}(x)$ as continuous functions on the
completion $\F_{\pi}[[\pi]]$.  
Wagner \cite[\S 5]{wagner} showed that 
the $r$-digit extension ${\mathcal E}_{i,T,r}(x)$ 
of the polynomials $E_{j,T,r}(x)$  
forms an orthonormal basis for all continuous functions 
from $\F_{\pi}[[\pi]]$ to its quotient field.  As a corollary 
Wagner showed the polynomials
$E_{j,T,r}(x)$ form an orthonormal basis for the 
$\F_r$-linear continuous functions on $\F_{\pi}[[\pi]]$.
We will prove these results in the reverse order, which 
seems more natural.

First we need a well-known lemma 
which is analogous to the mod $p^n$ periodicity of the 
binomial polynomials $\binom{x}{i} \colon 
\ZZ \rightarrow \ZZ/p\ZZ$ when $i < p^n$.

\begin{lemma}
\label{prep}
Let $\pi$ be irreducible in $\F_r[\UU]$, of degree $d$.  If 
$j < dn$, then $E_{j,\UU,r}(\pi^n g) \equiv 0 \bmod \pi$ for all 
$g$ in $\F_r[\UU]$. 
\end{lemma}

\begin{proof}
We may suppose $g \not= 0$, and have to show 
$\ord_{\pi}(e_j(\pi^n g)) > \ord_{\pi}(D_j)$, 
where $e_j$ and $D_j$ are the appropriate Carlitz objects 
on the ring $\F_r[\UU]$.

For integers $k \geq 0$, let $k \equiv R_k \bmod d$, where 
$0 \leq R_k \leq d-1$.  In particular, 
write $j = dQ + R_j$.  Since $\ord_{\pi}(D_k) = (r^k-r^{R_k})/(r^d-1)$ 
and $\deg(\pi^ng) > j$, 
\begin{eqnarray*}
\ord_\pi(e_j(f^ng)) & = & n + \ord_\pi(g) + \sum_{k=0}^{j-1}
(r-1)\ord_{f}(D_k) \\
& = & n + \ord_\pi(g) + \sum_{k=0}^{j-1}
(r-1)\left(\frac{r^k-r^{R_k}}{r^d-1}\right) \\
& = & n + \ord_\pi(g) + \frac{r^j-r^{R_j}}{r^d-1} - Q \\
& = & n + \ord_{\pi}(g) + \ord_{\pi}(D_j) - Q.
\end{eqnarray*}
Since $n > j/d \geq Q$, we're done.
\end{proof}

\begin{lemma}\label{1st}
Let $\pi$ be irreducible in $\F_r[\UU]$. 
The polynomials $E_{j,T,r}(x)$, 
viewed as $\F_r$-linear continuous functions on the completion 
$\widehat{\F_r[\UU]}_{\pi} = \F_{\pi}[[\pi]]$, are 
an orthonormal basis for all the $\F_r$-linear 
continuous maps from 
$\F_{\pi}[[\pi]]$ to $\F_{\pi}((\pi))$.
\end{lemma}

\begin{proof}
Let $d$ be the degree of $\pi$ and $n$ be any positive integer.
By Lemma \ref{prep}, for $j < dn$ 
$\overline{E}_{j,T,r}$ is a well-defined 
map from $\F_{r}[\UU]/(\pi^n)$ to 
$\F_{r}[\UU]/(\pi) \cong \F_{\pi}$.

For $0 \leq j, k \leq dn-1$, 
the $dn \times dn$ matrix $(E_{j,\UU,r}(\UU^k))$ 
is triangular with all diagonal entries equal to 1. 
Since $1, T, \dots, T^{dn-1}$ are an $\F_{r}$-basis of 
$\F_{\pi}[[\pi]]/(\pi^n) \cong \F_{r}[T]/(\pi^m)$, it follows that 
the $dn$ reduced functions 
$\overline{E}_{j,\UU,r}$ are an $\F_{\pi}$-basis of 
$\Hom_{\F_r}(\F_{\pi}[[\pi]]/(\pi^n),\F_{\pi})$.
Therefore $\{E_{j,T,r}\}_{j \geq 0}$ is an orthonormal 
basis of $\LC_{\F_r}(\F_{\pi}[[\pi]],\F_{\pi}((\pi)))$, as 
we wanted to show.
\end{proof}

\begin{theorem}
\label{Ejrcomp}
The Carlitz polynomials ${\mathcal E}_{i,T,r}$ in $\F_r[T]$ 
form an orthonormal basis for the continuous functions 
from $\F_{\pi}[[\pi]]$ to $\F_{\pi}((\pi))$ 
when $\pi$ is any irreducible in $\F_r[\UU]$.
\end{theorem}

\begin{proof}
To simplify the notation, 
we write $E_j$ for $E_{j,\UU,r}$.

Let $d$ be the degree of $\pi$.
Let $H = \cap_{j=0}^{dn-1} \Ker(\overline{E}_{j})$, so 
$(\pi^n) \subset H$ by Lemma \ref{prep}.  
We want to show $H = (\pi^n)$, 
and then we'll be done by Theorem \ref{rrq}.

By Lemma \ref{1st}, $\overline{E}_{0},\dots,\overline{E}_{dn-1}$ 
form an $\F_{\pi}$-basis of the $\F_r$-linear 
maps from $\F_r[\UU]/(\pi^n)$ to $\F_r[\UU]/(\pi)$. 
By an argument as in the proof of Theorem \ref{rrq}, 
this implies the functions 
$\overline{E}_{0},\dots,\overline{E}_{dn-1}$ separate 
the points of $\F_{r}[\UU]/(\pi^n)$.

So any element of $\F_r[\UU]$ which is killed by all 
$\overline{E}_{j}$ for $j \leq dn-1$ must be in 
$(\pi^n)$.  Therefore $H \subset (\pi^n)$. 
\end{proof}

The conclusion of Theorem \ref{Ejrcomp} 
is analogous to the role of the binomial 
polynomials $\binom{x}{n}$, which 
are an orthonormal basis of $C(\Zp,\Qp)$ for all primes $p$.
Note that 
the Carlitz polynomials in $\F_r[T]$ do not give an 
orthonormal basis in the completion at $1/T$, 
as they do not even take integral elements to 
integral elements at this place.  (Though see 
Car \cite{car} for a use of these polynomials 
$1/T$-adically to express a very large class of {\it entire}
power series. This class, for instance, includes the $L$-series of
Drinfeld modules etc.)

The next application of the digit principle 
(specifically, the proof of Theorem \ref{hassexp})
is the original motivation for this paper.

\begin{lemma}
\label{hyperbasis}
Let $K = \FF((\TT))$, $\OO = \FF[[\TT]]$.
The hyperdifferential functions 
$\{\DD_j\}_{j \geq 0}$ on $K$ are an orthonormal 
basis of $\LC_{\FF}(\OO,K)$.  
\end{lemma}

Here $\DD_{j} = \DD_{j,T}$.  We refer to these operators 
as functions in the lemma because we will later be considering 
their product in the sense of functions, not (via composites) in 
the sense of operators. 

Lemma \ref{hyperbasis} is independently due to
Jeong \cite{jeong} and 
 Snyder \cite{snyder}, with proofs different from the one 
we now give.

\begin{proof}
Composing the function $\DD_j$ with reduction mod $\TT$, we get an 
$\FF$-linear map $\overline{\DD}_j \colon \OO \rightarrow \FF$ 
whose kernel consists of power series with $T^j$-coefficient 0.
The reductions $\overline{\DD}_0,\dots,\overline{\DD}_{n-1}$ 
are well-defined elements of 
the $\FF$-dual space $(\FF[\TT]/\TT^n)^{*}$, 
and in fact are the dual basis 
to $1, \TT, \dots, \TT^{n-1}$.  We are done by Lemma \ref{orc}.
\end{proof}

\hfill

\noindent
{\bf Example.}
Let $\Phi_q$ be the $q$th power Frobenius on $\FF[[T]]$, 
so $\Phi_q(x) = x^q$.  
Then $\Phi_{q} = \sum_{j \geq 0} b_j \DD_{j}$ for some 
sequence $b_j$ in $\FF[[T]]$ tending to 0.   
Applying the binomial theorem to $T^{qn} = (T^q - T + T)^n$, 
we obtain $b_j = (T^q - T)^j$.  This expansion formula 
was already noted by Voloch \cite{vol}.

\hfill

An alternate proof of 
Lemma \ref{hyperbasis} 
comes from 
Lemma \ref{cbasis} and the observation 
of Jeong \cite{jeong2} that $\overline{\DD}_j = 
\overline{E}_j$ 
for all $j$.  The equality of these reduced functions 
contrasts with the rather different behavior of 
$\DD_j$ and $E_j$ as (linear) maps from $\FF[[T]]$ to $\FF[[T]]$: 
$\DD_j$ has an infinite-dimensional kernel and, 
as noted by Voloch \cite{vol}, $\DD_j$ is 
nowhere differentiable.

For an integer $i \geq 0$, let 
$$
i = c_0 + c_1q + \dots + c_{n-1}q^{n-1}
$$
be its base $q$ expansion, where $0 \leq c_j \leq q-1$.  Define
$$
\RD_i \eqdef \DD_{0}^{c_0}\DD_{1}^{c_1}\dots\DD_{n-1}^{c_{n-1}}, 
$$
where the product on the right is a product of continuous functions on 
$K = \FF((\TT))$, {\it not} a composite of operators.   
Avoiding this confusion is the reason we call the 
$\DD_j$ hyperdifferential functions, and not hyperdifferential operators, 
when they are viewed as functions.
Note $\DD_{j} = \RD_{q^j}$.

\begin{theorem}\label{hassexp}
Let $K = \FF((\TT)), \OO = \FF[[\TT]]$.
The sequence 
$\{\RD_i\}_{i \geq 0}$ is an orthonormal basis of $C(\OO,K)$.
\end{theorem}

\begin{proof}
Use Lemma \ref{hyperbasis} and the digit principle.
\end{proof}

In analogy to Theorem \ref{Ejrcomp}, 
we can use the (global) higher $\F_r$-derivation $\{\DD_{j,\UU}\}$ 
on $\F_r[\UU]$ to give an orthonormal basis for the 
continuous functions on completions of $\F_r[\UU]$.
While any completion of 
$\F_r(\UU)$ is a Laurent series 
field $\FF((u))$, 
the extension of the (global) hyperdifferential 
functions on $\F_r[\UU]$ to this completion will 
generally not be the hyperdifferential functions $\DD_{j,u}$ 
on $\FF((u))$ that are used in Lemma \ref{hyperbasis}.

The following result answers a question of Goss.

\begin{theorem}
\label{commonhyp}
Let $\pi$ be irreducible in $\F_r[T]$, of degree $d$, 
$\OO = \F_{\pi}[[\pi]]$ the corresponding completion 
at $\pi$, $K$ its fraction field.
The hyperdifferential functions 
$\DD_{j,\UU}$ on $\F_r[\UU]$, extended by continuity to 
$\OO$, give an 
orthonormal basis for the $\F_{\pi}$-linear 
functions from $\OO$ to $K$.
The extension of the $\DD_{j,T}$ by $r^d$-digit expansions 
gives  
an orthonormal basis of $C(\OO,K)$.
\end{theorem}

Note the digit extension of the sequence 
$\{\DD_{j,\UU}\}_{j \geq 0}$ to an orthonormal basis of $C(\OO,K)$ 
depends on the possible change in the residue field 
under completion, unlike the digit extension used in Theorem \ref{Ejrcomp}.

\begin{proof}
First we see why completion at $1/T$ is not being considered.
Write $S = 1/\UU$.  Since $\DD_{j,\UU}(S^m) = \binom{-m}{j}S^{m+j}$, 
$\DD_{j,\UU}$ has image in $S^{j}\F_r[[V]]$.  So these functions 
are not orthonormal on the completion at $1/T$.

Now we look at the completion $\OO = \widehat{\F_r[T]}_{\pi}$.
To establish the first claim of the theorem,
it suffices by the digit principle to check 
the $\DD_{j,\UU}$ are an orthonormal 
basis of 
$\LC_{\F_{\pi}}(\OO,K)$.
We've already checked in Theorem \ref{Djrcomp} that they belong to this space.

By (\ref{hypercong2}) and continuity, 
the reduced functions 
$\overline{\DD}_{j,\UU} \colon \OO \rightarrow \F_{\pi} \cong 
\F_{r}[\UU]/(\pi)$, 
for $0 \leq j \leq n-1$, annihilate the ideal $(\pi^n)$.
We now check the corresponding functions on $\OO/(\pi^n)$ 
are a basis of the $\F_{\pi}$-dual space, which will end the proof 
by Lemma \ref{orc}.

Consider the effect of these $n$ functions on the basis $1,\pi,
\dots,\pi^{n-1}$. 
Since $\DD_{j,\UU}(\pi^n) \equiv 0 \bmod \pi$ for $j < n$, 
the $n \times n$ matrix $(\overline{\DD}_{j,\UU}(\pi^k))$ is triangular. 
Since $\pi'(T) \not\equiv 0 \bmod \pi$, (\ref{hypercong1}) shows 
$\DD_{j,\UU}(\pi^j) \not\equiv 0 \bmod \pi$, so the diagonal 
entries are all nonzero.
So this matrix is invertible.
\end{proof}

For all (monic) irreducibles of a fixed degree in $\F_r[\UU]$, 
Theorem \ref{commonhyp}
gives a single family of 
non-polynomial functions which serves as   
an orthonormal basis of the space of continuous functions on the 
completion at each of these irreducibles.

For a monic irreducible $\pi$ in $\F_r[T]$, we'd like a 
Chain Rule formula 
for computing the effect of all the $\DD_{j,T}$ on the completion 
$\F_{\pi}[[\pi]]$ in terms of both  
the effect of all the $\DD_{j,\pi}$ and 
the data $\DD_{j,T}(\pi)$.

It suffices 
to give a formula for $\DD_{j,\UU}(\pi^n)$ when $j \geq 1$, 
which follows from the 
Leibniz rule (\ref{mprod}) with multiple factors:
\begin{equation}\label{teichcr}
\DD_{j,\UU}(\pi^n) = \sum_{\Atop{k_1 + \dots + k_n = j}{k_1,\dots,k_n \geq 0}}
\DD_{k_1,\UU}(\pi)\cdots \DD_{k_n,\UU}(\pi) = 
\sum_{i=1}^{j} \binom{n}{i}\pi^{n-i} \sum_{\Atop{k_1 + \dots + k_i = j}
{k_1,\dots,k_i \geq 1}} \DD_{k_1,\UU}(\pi) \cdots \DD_{k_i,\UU}(\pi).
\end{equation}
The second sum on the right side 
simply collects together all tuples $(k_1,\dots,k_n)$ from the first 
sum having 
the same number $i$ of positive coordinates.  The remaining coordinates
in the tuple 
are 0, and this contributes a factor of $\pi^{n-i}$.

Extending (\ref{teichcr}) by linearity and continuity
gives a direct Chain Rule for $\DD_{j,\UU}(f(\pi))$ for any 
$f(\pi) \in \F_{\pi}[[\pi]]$:
$$
\DD_{j,\UU}(f(\pi)) = 
\sum_{i=1}^{j} \DD_{i,\pi}(f(\pi))\sum_{\Atop{k_1 + \dots + k_i = j}
{k_1,\dots,k_i \geq 1}} \DD_{k_1,\UU}(\pi)\cdots\DD_{k_i,\UU}(\pi).
$$
This is due to Teichm\"uller \cite[Equation 6]{teichm}.

For a local field $K$ of positive characteristic, 
the digit principle provides us with clearer picture of 
how generally linear functions in $C(\OO,K)$ can be built up to an 
orthonormal basis, and also suggests alternate 
``canonical'' isomorphisms between nonarchimedean 
measures and formal divided power series. 
A correspondence between measures and such series arises 
because of the addition formula for Carlitz polynomials:
\begin{equation}\label{addform}
{\mathcal E}_{i}(x+y) = \sum_{j+k=i} \binom{i}{j}
{\mathcal E}_{j}(x){\mathcal E}_{k}(y).
\end{equation}
This formula motivates the assignment to a measure $\nu$ on 
$\FF[[\TT]]$ the formal divided power series $\sum_{i \geq 0} 
(\int_{\FF[[\TT]]} {\mathcal E}_i(x)\,d\nu) (X^i/i!)$.  
A useful property of 
this correspondence is that
convolution of measures corresponds to the simpler operation of 
multiplication of the corresponding 
series.  (This is analogous to the effect of the 
Fourier transform, which converts convolution into 
multiplication.)
Since the addition formula for ${\mathcal E}_i(x+y)$ 
follows purely from the construction of the
Carlitz polynomials ${\mathcal E}_i$ in 
terms of $\FF$-linear functions and digit expansions (cf. 
Goss \cite[Prop. 3.2.1]{goss}), we can replace the Carlitz basis 
with other orthonormal bases in characteristic $p$ which are 
constructed by the digit principle.  
Namely, 
if $\{e_j\}$ is {\it any} fixed orthonormal 
basis of $\LC_{\FF}(\FF[[\TT]],\FF((\TT)))$ and $\{f_i\}$ is 
the orthonormal basis of $C(\FF[[\TT]],\FF((\TT)))$ 
constructed from the $e_j$ by $q$-digits, then 
attaching to an $\FF((\TT))$-valued measure $\nu$ the 
formal divided power series 
$\sum_{i \geq 0} (\int_{\FF[[\TT]]} f_i(x)\,d\nu) (X^i/i!)$ 
converts convolution of measures into products of series.
(If $q = r^d$, this also applies to $r$-digit extensions of an orthonormal 
basis of $\LC_{\F_r}(\FF[[T]],\FF((T)))$ which satisfies 
the kernel hypothesis of Theorem \ref{rrq}.)

We now turn to some applications of the digit principle 
in characteristic 0, in the form 
of Theorem \ref{char0}.

\begin{theorem}\label{coll}
For $m \geq 0$, write $m = c_0 + c_1p + \dots + c_kp^k$ where 
$0 \leq c_j \leq p-1$.  Set 
$$
\left\{\Atop{x}{m}\right\} \eqdef
\binom{x}{1}^{c_0}\binom{x}{p}^{c_1}\cdots \binom{x}{p^{k}}^{c_{k}}.
$$
The functions $\left\{\Atop{x}{m}\right\}$ are an orthonormal 
basis of $C(\Zp,\Qp)$.
\end{theorem}

\begin{proof}
For $0 \leq i \leq p^n-1$ and $x, y \in \Zp$, 
$$
x \equiv y \bmod p^n \Longrightarrow 
(1 + T)^x \equiv (1 + T)^y \bmod (p,T^{p^n}) \Longrightarrow 
\binom{x}{i} \equiv \binom{y}{i} 
\bmod p, 
$$
so the $p^n$ functions $\binom{x}{i}$ are well-defined 
maps from $\ZZ/p^n\ZZ$ to $\ZZ/p\ZZ$. 
To prove the theorem, it suffices by Theorem \ref{char0} to show that 
for each $x \in \ZZ/p^n\ZZ$, the sequence
$$
\binom{x}{1} \bmod p, \binom{x}{p} \bmod p, \dots, \binom{x}{p^{n-1}} \bmod p
$$
determines $x$.  Writing $x \equiv d_0 + d_1p + \dots + d_{n-1}p^{n-1}$ 
with $0 \leq d_j \leq p-1$, Lucas' congruence implies 
$\binom{x}{p^j} \equiv d_j \bmod p$, so we're done.
\end{proof}

The orthonormal basis $\left\{\Atop{x}{m}\right\}$ of $C(\Zp,\Qp)$ 
is similar in appearance to the 
Carlitz basis for $C(\FF[[\TT]],\FF((\TT)))$, 
but it does not have algebraic features as 
convenient in characteristic 0 
as the usual Mahler basis $\binom{x}{n}$ of $C(\Zp,\Qp)$.

Mahler's basic theorem about the binomial coefficient 
functions is a consequence of the previous theorem, as follows.

\begin{corollary}
\label{mpf}
The functions $\binom{x}{n}$ are an orthonormal 
basis of $C(\Zp,\Qp)$.
\end{corollary}

\begin{proof}
Since each $\left\{\Atop{x}{n}\right\}$ has degree $n$ and 
$\binom{x}{n}$ sends $\Zp$ to $\Zp$, the transition matrix 
from $\{\Atop{x}{0}\},\dots,\{\Atop{x}{n}\}$ to 
$\binom{x}{0},\dots,\binom{x}{n}$ is triangular over $\Zp$ with 
diagonal entries
$$
\frac{i!}{(1!)^{c_0}(p!)^{c_1}\cdots (p^k)!^{c_k}},
$$
where $i = c_0 + c_1p + \dots + c_kp^k$, $0 \leq c_j \leq p-1$.
This ratio is a $p$-adic unit, so the reduced functions 
$\binom{x}{n} \bmod p$ 
are a basis of $C(\Zp,\Fp)$.  We are done by Lemma \ref{orc}.
\end{proof}

Writing $x = d_0 + d_1p + d_2p^2 + \dots$, with $0 \leq d_j \leq p-1$, 
we can compare the reductions of $\binom{x}{m}$ and $\left\{\Atop{x}{m}
\right\}$ as functions from $\Zp$ to $\Fp$:
$$
\binom{x}{m} \equiv 
\binom{d_0}{c_0} \cdots \binom{d_k}{c_k} \bmod p,  \ \ 
\left\{\Atop{x}{m}\right\} \equiv d_0^{c_0} \cdots d_k^{c_k} \bmod p.
$$

In light of the diagonal matrix entries in the proof of Corollary 
\ref{mpf}, 
probably the closest analogue for hyperdifferential operators 
on $\Fp[[T]]$ is
$$
\DD_{1}^{\circ c_0} \circ 
\DD_{p}^{\circ c_1} \circ \cdots \circ \DD_{p^k}^{\circ c_k} = 
\frac{i!}{1!^{c_0}(p!)^{c_1}\cdots (p^k!)^{c_k}}
\DD_{i},
$$
where $i = c_0 + c_1p + \dots + c_kp^k$, $0 \leq c_j \leq p-1$. 
Here $\DD_{j}^{\circ c}$ is the $c$-fold composite of 
$\DD_{j}$.  So this only provides us with another 
basis for the linear continuous functions.

The orthonormal bases $\binom{x}{n}$ and $\left\{\Atop{x}{n}\right\}$ 
consist of 
polynomials.  A criterion for a sequence of polynomials $P_n(x)$ to be an 
orthonormal basis of $C(\OO,K)$ can be given in 
terms of degrees and leading coefficients, avoiding the 
appeal to Lemma \ref{orc} which we consistently make.
See Cahen and Chabert \cite{cc}, De Smedt \cite{desmedt}, or
Tateyama \cite{tateyama}.
As shown by Yang \cite{yang}, the conditions 
for analyticity or local analyticity for functions in 
$C(\Zp,\Qp)$, in terms of Mahler 
coefficients as given by 
Amice \cite{amice}, carry over to these 
polynomial orthonormal bases $P_n(x)$.

The construction in Theorem \ref{coll} 
is formulated more generally by Tateyama \cite{tateyama} 
using coefficient functions arising 
from Lubin-Tate formal groups.  Let $K$ be a local field, 
with ring of integers $\OO$ and residue field size $q$. 
Fix a uniformizer $\pi$ and 
a Lubin-Tate formal group 
$F/\OO$ associated to some 
Frobenius power series $[\pi](X) \in \OO[[X]]$.  
We write $[a](X) = [a]_{F}(X)$ for the endomorphism of $F$ attached 
to each $a \in \OO$.  Write 
\begin{equation}\label{tdef}
[a](X) = \sum_{n \geq 1} C_{n,F}(a)X^n,
\end{equation}
which defines functions $C_{n,F} \colon \OO 
\rightarrow \OO$.  

In characteristic 0, letting $\lambda_F \colon F \rightarrow {\mathbb G}_a$ 
be the unique normalized logarithm, with $\exp_F$ its composition 
inverse, the equation $\lambda_{F}([a](X)) = a\lambda_F(X)$ 
leads to $[a](X) = \exp_{F}(a\lambda_{F}(X))$.  
Comparing with (\ref{tdef}) shows $C_{n,F}(a)$ is a 
polynomial function of $a$, with degree at most $n$. 
(Using the formal group law and (\ref{tdef}) alone, one 
could check $C_{n,F}(a)$ is continuous in $a$, but it's 
easier to obtain this from knowing that $C_{n,F}$ is 
actually a polynomial.)  
Tateyama observes that while there is no 
unique normalized logarithm for $F$ 
in characteristic 
$p$, in all characteristics we can 
carry out the same argument from characteristic 0 by using 
Wiles' construction of a 
logarithm, given by the coefficient-wise limit formula
$$
\lambda_{F}(X) \eqdef \lim_{n \rightarrow \infty} \frac{[\pi^n](X)}{\pi^n} = 
X + \dots.
$$
By an explicit check, this particular logarithm 
satisfies $\lambda_F([a](X)) = a\lambda_F(X)$ 
even in characteristic $p$, so 
$C_{n,F}(a)$ is a polynomial 
in $a$ (of degree at most $n$) in all cases.

\hfill

\noindent
{\bf Example.}
$F/\ZZ_2$ is the Lubin-Tate group attached to $[2](X) = X^2 + 2X = 
(1+X)^2 - 1$, so $F = {\mathbb G}_m$ and $C_{n,F}(a) = \binom{a}{n}$.

\hfill

\hfill

\noindent
{\bf Example.}
$F$ is the Lubin-Tate group over $\FF[[T]]$ attached to the 
series $[T](X) = X^q + TX$, i.e., $F$ is the Carlitz module. 
Then 
$$
C_{n}(a) = 
\begin{cases}
E_{k}(a), & \text{ if } n = q^k; \\
0, & \text{ if } n \text{ is not a power of } q.
\end{cases}
$$

\hfill

\begin{theorem}
\label{tatey}
Let $\OO$ be the integer ring of a local field, $\F$ the residue field, $q 
= \#\FF$.  For a Lubin-Tate group $F/\OO$,
the polynomials  
$$
C_{1,F}(x)^{c_0}C_{q,F}(x)^{c_1}\cdots C_{q^{k-1},F}(x)^{c_{k-1}},  \ \ \ 
k \geq 1,  \ \ 0 \leq c_j \leq q-1,
$$
form an orthonormal basis of $C(\OO,K)$.
\end{theorem}

While Tateyama \cite{tateyama} proves this by 
a criterion on polynomial orthonormal bases, we'll use 
the digit principle instead.  Both 
Theorems \ref{cw1} and \ref{coll} are special cases.

\begin{proof}
Let $[\pi](X)$ be the Frobenius series attached to $F$. 
Since $[\pi](X) \equiv X^q \bmod \pi$, 
$$
[\pi^{j+1}a](X) \equiv ([a](X))^{q^{j+1}} \bmod \pi \equiv 
0 \bmod (\pi, X^{q^{j+1}})
$$ 
for 
all $a \in \OO$.  So 
for $m < q^{j+1}$, $\overline{C}_{m} \colon \OO \rightarrow \F$ 
annihilates $(\pi^{j+1})$.  Taking $m = 1, q, \dots, q^{n-1}$, 
we will show the induced map
\begin{equation}\label{tymap}
\OO/\pi^n \rightarrow \F^n \ \ 
\text{ given by } \ \ x \mapsto (\overline{C}_1(x),\overline{C}_{q}(x),\dots,
\overline{C}_{q^{n-1}}(x))
\end{equation}
is a bijection, so we'd done by the digit principle.

For $a \in \OO$, 
$[\pi^{j}a] \equiv ([a](X))^{q^j} \bmod \pi$, so 
\begin{equation}\label{flucas}
C_{q^j}(\pi^{j}a) \equiv a^{q^j} \equiv a \bmod \pi.
\end{equation}
Therefore $C_{q^j}$ recovers the $\pi^j$-coefficient 
of any element of the ideal $(\pi^j)$.  Take $j = 0, 1, \dots, n-1$ 
successively in the congruence 
(\ref{flucas}), which is a formal group version of 
a weak form of Lucas' congruence: $\binom{dp^j}{p^j} \equiv d \bmod p$ 
for $0 \leq d \leq p-1$. 
So we see that (\ref{tymap}) is a bijection.
\end{proof}

Our next example is an orthonormal basis due to Baker \cite{baker}, 
consisting
not of polynomials, 
but of locally constant functions 
taking values in the Teichm\"uller representatives.

\begin{theorem}\label{bbasis}
Let $K$ be any local field, $\OO$ its ring of integers, $\pi$ a 
fixed uniformizer of $\OO$, $q = \#\OO/\pi$. For each $x \in \OO$, write 
$$
x = \sum_{j \geq 0} \omega_j(x)\pi^j, 
$$
where $\omega_{j}(x)$ is a Teichm\"uller representative.  

For $m \geq 0$ with $m = c_0 + c_1q + \dots + c_kq^k$, $0 \leq c_j \leq q-1$, 
let 
$$
\Omeg_m(x) \eqdef 
\omega_{0}(x)^{c_0} \omega_{1}(x)^{c_1} \cdots \omega_{k}(x)^{c_{k}}.
$$
The functions $\Omeg_m(x)$ for $m \geq 0$ are an orthonormal basis 
of $C(\OO,K)$.
\end{theorem}

\begin{proof}
The functions $\omega_i(x)$ for 
$i = 0,\dots,n-1$ obviously separate the points 
of $\OO/\pi^n$.  Now apply Theorem \ref{char0}.
\end{proof}

Baker's proof of Theorem \ref{bbasis} differs from ours in the demonstration 
that the functions $\overline{\Omeg}_{0}(x),
\overline{\Omeg}_1(x),\dots,
\overline{\Omeg}_{q^n-1}(x)$ (for each $n$) are linearly 
independent in $\Maps(\OO/\pi^n,\FF)$.  While the argument 
given here using the digit 
principle shows these functions (in a $q^n$-dimensional 
space) are linearly independent because they are 
a spanning set,
Baker shows the linear independence by a technical direct 
calculation.  (He also provides a set of 
polynomial functions on $\OO$ whose reductions 
coincide with the $\overline{\omega}_i$.)

Referring to the sequence $\Omeg_m(x)$ as the Teichm\"uller 
basis may create confusion with the expansion 
of elements of $K$ in terms of Teichm\"uller representatives, so we 
call the sequence $\Omeg_m(x)$ the Baker basis.
Note $\Omeg_{q^k}(x) = \omega_k(x)$. 
For an integer $R \geq 0$, 
Baker \cite{baker} writes $\Omeg_R(x)$ as $\omega^{R}(x)$.

For $p$ odd and $K = \Qp$, $\Omeg_{(p-1)/2}(x) = \omega_0(x)^{(p-1)/2} = 
(\frac{x}{p})$ is the Legendre 
symbol.  Higher power residue symbols are in the Baker basis 
for suitable finite extensions of $\Qp$.

\section{The Baker Basis and the Tate Algebra}

We continue with the same notation as at the end of 
Section \ref{apps}.  
In particular, $K$ is any local field and $\F$ is 
its residue field, with $q$ the size of $\F$.
Since $\Omeg_m(x)$ depends on the choice of $\pi$ for $m \geq q$, 
a better notation is $\Omeg_{m,\pi}(x)$. Since the functions 
$\omega_j(x) = \Omeg_{q^j}(x)$ 
depend on $\pi$ (except when $j = 0$), we could 
write them as $\omega_{j,\pi}(x)$.

As an example of an expansion in the Baker basis,
the expansion of each $x \in \OO$ 
using Teichm\"uller representatives 
amounts to giving the Baker expansion of the identity function:
\begin{equation}\label{xexp}
x = \sum_{j \geq 0} \omega_{j}(x)\pi^j \Longrightarrow 
x = \sum_{j \geq 0} \pi^j\Omeg_{q^j}(x).
\end{equation}
Therefore
$$
x^2 = \sum_{i,j \geq 0} \pi^{i+j}\omega_{i}(x)\omega_j(x) = 
\sum_{i, j \geq 0} \pi^{i+j}\Omeg_{q^i + q^j}(x), 
$$
which is a Baker expansion except if $q = 2$. 
In that case we simplify with the rule
$\omega_i(x)\omega_i(x) = \omega_i(x)$.

Because 
the functions $\Omeg_m(x)$ behave very simply under multiplication, 
we'll see below (Theorem \ref{tate}) that 
they elucidate the 
structure of $C(\OO,K)$ as a $K$-Banach algebra, in terms of 
the ``infinite-dimensional'' Tate algebra 
$T_{\infty}(K) \eqdef K\langle X_1, X_2,\dots\rangle$. 
As a set, $T_{\infty}(K)$ consists of the formal power series  
$f(\underline{X}) = \sum_{i} a_{i}\underline{X}^i \in K[[X_1,X_2,\dots]]$ in 
countably many indeterminates such that $a_i \rightarrow 0$ 
as $i \rightarrow \infty$.  That is, 
for any $\varepsilon > 0$, $|a_i| < \varepsilon$ 
for all but finitely many $i$.  
(The indices $i$ run through sequences in 
$\NN^{(\infty)} = \oplus_{n \geq 0} \NN$, and 
$\underline{X}^i$ denotes a monomial of several 
variables, such as $X_1^{i_1}\cdots X_m^{i_m}$.)
Note that $\sum_{j \geq 1} X_j$ is 
not in $T_{\infty}(K)$.  

The set $T_{\infty}(K)$ has a 
natural $K$-algebra structure.
We topologize $T_{\infty}(K)$ 
using the sup-norm on coefficients,
$$
|f(\underline{X})| \eqdef \sup_{i} |a_i|.
$$
So the unit ball of $T_{\infty}(K)$ 
is the $\mm(X_1,X_2,\dots)$-adic completion of 
the polynomial algebra $\OO[X_1,X_2,\dots]$, 
where $\mm$ is the maximal ideal of $\OO$.

The algebra $T_{\infty}(K)$ shares some properties with the 
more traditional finite-dimensional 
Tate algebras $T_{n}(K) = K\langle X_1,\dots,X_n\rangle$,
e.g., there are R\"uckert Division and Weierstrass Preparation 
Theorems for $T_{\infty}(K)$, from which one can show 
$T_{\infty}(K)$ has  
unique factorization (but not by induction on 
the number of variables, as is traditionally the case for  
$T_{n}(K)$).  
This will not be needed for what follows, so we defer the proof to 
a later paper.

There are differences between $T_{\infty}(K)$ and 
$T_{n}(K)$, the most 
obvious being that $T_{\infty}(K)$ is 
not noetherian. The 
ideal $(X_1,X_2,X_3,\dots)$ is not finitely generated, 
and also not closed, e.g., 
the series $\sum \pi^jX_j$ is in the closure of the 
ideal but not in the ideal.  It is easy to write down many 
more non-closed ideals of $T_{\infty}(K)$ in a similar manner.

For any closed ideal $I$ of $T_{\infty}(K)$, we equip 
$T_{\infty}(K)/I$ with the residue norm:
$|f \bmod I|_{\res} = \inf |f + h|$, where the infimum is taken 
as $h$ runs over 
$I$.  This residue norm makes $T_{\infty}(K)/I$ a $K$-Banach algebra 
whose norm topology is the quotient topology
\cite[Prop. 1.1.6/1, 1.1.7/3]{bgr}.

Recall the residue field $\F$ of $K$ has size $q$. 
Call a series $\sum a_i\underline{X}^i \in T_{\infty}(K)$ 
{\it q-simplified} if the exponents in every nonzero monomial term
are all at most $q-1$.  Let $I_q(K)$ be the closure 
of the ideal generated by all $X_j^q-X_j$.

\begin{lemma}\label{norm}
Every congruence class in $T_{\infty}(K)/I_q(K)$ has a unique 
q-simplified representative, and if $f \equiv g \bmod I_q(K)$ 
with $g$ being $q$-simplified, then $|f \bmod I_q(K)|_{\res} = |g|$.
\end{lemma}

\begin{proof}
For $m \geq 0$ and $Y$ an indeterminate, we can write 
in $\ZZ[Y]$
$$
Y^m \equiv Y^{m'} \bmod (Y^q - Y)
$$
for some (unique) $m' \leq q-1$.  So for any monomial 
$\underline{X}^i$ in the variables $X_1,\dots,X_n$, we can write
$$
\underline{X}^i = \underline{X}^{i'} + h_i(\underline{X}),
$$
where all the exponents in $i'$ are $\leq q-1$ and 
$h_i \in \ZZ[X_1,\dots,X_n]$ 
is in the ideal generated by $X_1^q - X_1,\dots,X_n^q-X_n$.  
Viewing this 
equation in $K[X_1,\dots,X_n]$, note $|h_i| \leq 1$. 
So for any series $f = \sum a_i\underline{X}^i \in 
T_{\infty}(K)$, we can write $f = g + h$ where $g$ is $q$-simplified 
and $h \in I_q(K)$.  By construction, each coefficient 
of $g$ is a (convergent) 
sum of coefficients of $f$, so $|g| \leq |f|$.

We have proved existence of a $q$-simplified series in 
each class of $T_{\infty}(K)/I_q(K)$.  
Provided we show uniqueness, we then vary $f$ within a 
congruence class to see that $|f \bmod I_q(K)|_{\res} = |g|$.

For uniqueness, 
it suffices to show the only $q$-simplified series in $I_{q}(K)$ 
is 0.  Let $g \in I_{q}(K)$ be $q$-simplified and nonzero. 
Scaling, we may assume $|g| = 1$.  Then 
$\overline{g} = g \bmod \mm$ is a nonzero polynomial in 
$\F[X_1,X_2,\dots]$, say $\overline{g} \in 
\F[X_1,\dots,X_N]$.  Since $g \in I_{q}(K)$, 
$\overline{g}$ vanishes at all points in $\F^N$. 
Since $\#\F = q$ and all exponents of 
$\overline{g}$ are at most $q-1$, 
we must have $\overline{g} = 0$, which is a contradiction.
\end{proof}

To make a connection between $T_{\infty}(K)$ and $C(\OO,K)$, it 
is convenient to index the variables in $T_{\infty}(K)$ 
starting at 0, so we write $T_{\infty}(K) = K\langle X_0,X_1,\dots\rangle$.
The reason for this adjustment is 
that the first term in $\pi$-adic expansions in $\OO$ 
is naturally indexed by 0, not by 1.

\begin{theorem}\label{tate}
For any local field $K$, with ring of integers $\OO$ and residue field 
of size $q$, there is a $K$-Banach algebra isometric isomorphism
$$
K\langle X_0,X_1,\dots,\rangle/I_q(K) \cong C(\OO,K).
$$
This isomorphism depends on a choice of uniformizer of $\OO$.
\end{theorem}

\begin{proof}
Fix a uniformizer $\pi$ of $\OO$.  Using the basis $\Omeg_{m,\pi}(x)$, 
we can express any $f \in C(\OO,K)$ uniquely in the form
\begin{equation}\label{bexp}
f(x) = \sum_{m \geq 0} a_{m}\Omeg_{m,\pi}(x) = 
\sum_{k \geq 0} \sum_{c_0,\dots,c_k \leq q-1} 
a_{c_0+\dots+c_{k}q^k}
\omega_0(x)^{c_0}\cdots \omega_{k,\pi}(x)^{c_k}.
\end{equation}

The map $K[X_0,X_1,\dots] \rightarrow C(\OO,K)$ sending 
$X_n$ to $\omega_n(x)$ extends by 
continuity to a $K$-algebra homomorphism 
$T_{\infty}(K) \rightarrow C(\OO,K)$.
By (\ref{bexp}), this map is surjective.  Obviously 
each $X_j^q - X_j$ is in the kernel, so we get an 
induced surjection
$\psi \colon T_{\infty}(K)/I_q(K) \rightarrow C(\OO,K)$.  
Since the Baker basis of $C(\OO,K)$ is orthonormal, 
we focus our attention on $q$-simplified Tate series 
and see that 
$\psi$ is an isometry by Lemma \ref{norm}.
Therefore $\psi$ is an isometric isomorphism of 
$K$-Banach algebras.  As $\Omeg_{m,\pi}(x)$ depends on 
$\pi$, so does the isomorphism we've constructed.
\end{proof}

The proof shows there is a $K$-Banach space (but not $K$-Banach algebra) 
isomorphism between 
$C(\OO,K)$ and the space of $q$-simplified series in 
$T_{\infty}(K)$.  For that matter, any $K$-Banach algebra with 
a (countable) orthonormal basis as a $K$-Banach space will be 
algebraically a quotient of $T_{\infty}(K)$.  The special 
aspect of the above proof is that we can identify the 
corresponding ideal 
very simply and check the isomorphism is an isometry as well.

When $K = \Qp$, the Mahler basis 
$\binom{x}{n}$ suggests a picture of the 
algebra structure of $C(\Zp,\Qp)$ which is more complicated than 
what we see by Theorem \ref{tate}, since the functions $\binom{x}{n}$ 
satisfy the complicated multiplicative relations
$$
\binom{x}{i}\binom{x}{j} = \sum_{j \leq k \leq i+j} 
\binom{k}{i}\binom{i}{k-j}\binom{x}{k} = 
 \sum_{i \leq k \leq i+j} 
\binom{k}{j}\binom{j}{k-j}\binom{x}{k}.
$$

For examples of how some continuous functions 
look under the isomorphism of Theorem \ref{tate}, 
we simply have to remember that Theorem \ref{tate} identifies 
the function $\omega_{j,\pi}(x)$ with $X_j$.  
So the characteristic 
function of $\OO^{\times}$ corresponds to 
$X_0^{q-1}$ and the characteristic function 
of $\mm$ corresponds to 
$1 - X_0^{q-1}$.  As a check, the product of these functions 
in $T_{\infty}(K)/I_q(K)$ is
$$
X_{0}^{q-1}(1 - X_{0}^{q-1}) = X_{0}^{q-2}(X_0 - X_0^{q}) = 0, 
$$
as expected.  The characteristic function of 
the ball $a + \pi^n\OO$ corresponds to 
$$
\prod_{j=0}^{n-1} ( 1- (X_j - \omega_{j,\pi}(a))^{q-1}).
$$
In particular, the space of locally constant functions 
from $\OO$ to $K$ is $\oplus_{m \geq 0} K\Omeg_m$, 
which corresponds to the polynomial algebra in $T_{\infty}(K)/I_q(K)$ 
generated over $K$ by the $X_j$.

By (\ref{xexp}), the subset of $T_{\infty}(K)/I_q(K)$ 
corresponding to the $K$-analytic functions which converge on the closed 
unit disc in $K$ is the space of power series $\sum b_jY_{\pi}^j$, where 
$Y_{\pi} = \sum_{j \geq 0} \pi^jX_j \bmod I_{q}(K)$ and 
$b_j \rightarrow 0$.
This identification is not topological, since the usual topology on the 
space of $K$-analytic functions is not that coming from its 
embedding into the continuous functions.

The isomorphism in Theorem \ref{tate} is analogous to 
(\ref{symis}), and in fact recovers (\ref{symis}).
Namely, from Theorem \ref{tate} we obtain (with $\F$ the residue 
field of $K$)
$$
C(\OO,\F) \cong \F[X_0,X_1,\dots,]/(X_0^q-X_0,X_1^q-X_1,\dots), 
$$
from which it follows (keeping in mind the link between 
$X_j$ and $\omega_{j,\pi}(x)$) that
$$
\Maps(\OO/\mm^n,\F) = 
C(\OO/\mm^n,\F) \cong \F[X_0,\dots,X_{n-1}]/(X_0^q-X_0,\dots,
X_{n-1}^q - X_{n-1}), 
$$
which is essentially (\ref{symis}).

\begin{corollary}\label{pmax}
Every closed prime ideal of $C(\OO,K)$ is a maximal ideal of the form 
$M_x \eqdef \{f : f(x) = 0\}$, as $x$ varies over $\OO$.
\end{corollary}

\begin{proof}
Let $\mathfrak p$ be a closed prime ideal of $C(\OO,K)$, 
$q$ the size of the residue field of $K$.  Viewing 
$\mathfrak p$ as a prime ideal of $T_{\infty}(K)$ which contains 
$I_{q}(K)$, the containment $X_j^q - X_j \in \mathfrak p$ implies 
$X_j - \alpha_j \in \mathfrak p$ for a unique Teichm\"uller 
representative $\alpha_j$ of $K$.  Therefore $\mathfrak p$ contains 
the closure of $(X_1-\alpha_1, X_2 - \alpha_2,\dots)$, which 
is the maximal ideal $M_x$ for $x = \sum \alpha_j\pi^j$.
\end{proof}

Since all maximal ideals in a Banach algebra are closed, 
Corollary \ref{pmax} classifies all maximal ideals $M$ of 
$C(\OO,K)$, so we obtain 
$$
\sup_{x \in \OO} |f(x)| = \sup_{M} |f \bmod M|.
$$
Any $K$-algebra homomorphism from 
a $K$-Banach algebra to $C(\OO,K)$ is therefore continuous 
\cite[Prop. 3.8.2/3]{bgr}.  In particular, by the Open 
Mapping Theorem  
$C(\OO,K)$ has only one $K$-Banach algebra topology.

These calculations related to maximal ideals of $C(\OO,K)$ 
are special cases of what is known concerning 
$C(X,F)$ for any 
complete nonarchimedean field $F$ and any 
compact Hausdorff totally disconnected space $X$, 
e.g., all maximal ideals of $C(X,F)$ are of the form 
$\mm_x = \{f : f(x) = 0\}$.  See
van Rooij \cite[Chap. 6]{vr}, where it is also shown 
that for compact Hausdorff totally disconnected spaces $X$ and $Y$, 
there is a bijection between continuous functions from 
$X$ to $Y$ and $F$-algebra 
homomorphisms from $C(Y,F)$ to $C(X,F)$.

For any complete extension field $L$ of $K$, such as a 
completion of an algebraic closure of $K$, taking 
completed tensor products 
shows 
$$
C(\OO,L) \cong T_{\infty}(L)/I_q(L)
$$
as $L$-Banach algebras.
Here $\OO$ still denotes the integers of $K$.

For a fixed uniformizer $\pi$ of $K$, the $\pi$-adic expansion of 
elements of $\OO$ using Teichm\"uller representatives 
gives a homeomorphism $x \mapsto (\omega_{j,\pi}(x))_{j \geq 0}$ 
from $\OO$ to the product of 
countably many copies of the finite discrete space 
$\Teich(K) \eqdef \{z \in K : z^q = z\}$, with the product 
space having the product topology.  Let $K^a$ denote the 
algebraic closure of $K$.
We can think of $T_{\infty}(K)$ as the space of 
$K$-analytic functions on the infinite-dimensional 
unit 
ball $B^{\infty}(K^{a}) \eqdef \{(x_j) : x_j \in K^{a}, |x_j| \leq 1\}$, 
with the caveat that if not all coordinates 
$x_j$ of a point $x = (x_j)$ are in a common finite extension of 
$K$, then the value at $x$ of a series in $T_{\infty}(K)$ 
may need to be viewed in the completion $\widehat{K^a}$.
So we can think of $C(\OO,K)$, a space of continuous functions 
on $\OO$,  
roughly as the space of 
$K$-analytic functions on the subset of points $(x_j)$ in
$B^{\infty}(K^a)$
cut out by the equations $x_j^q = x_j$.  
The task of making this formulation more precise 
suggests trying to develop some 
type of ``infinite-dimensional'' 
rigid analysis using model spaces like  
$B^{\infty}(K^{a})$.

\hfill

\vfill

\end{document}